\documentclass[reqno,centertags, 12pt, draft]{article}

\usepackage{amsmath,amsthm,amscd,amssymb}

\usepackage{latexsym}

\usepackage{mathrsfs}

\newcommand{\bbR}{{\mathbb{R}}}

\newcommand{\bbZ}{{\mathbb{Z}}}

\newcommand{\bbC}{{\mathbb{C}}}

\newcommand{\scrV}{\mathscr{V}}

\newcommand{\no}{\nonumber}

\newcommand{\ti}{\tilde  }

\newcommand{\beq}{\begin{equation}}

\newcommand{\eeq}{\end{equation}}

\newcommand{\ba}{\begin{align}}

\newcommand{\ea}{\end{align}}

\newcommand{\vphi}{\varphi}

\allowdisplaybreaks

\numberwithin{equation}{section}

\newtheorem{theorem}{Theorem}[section]

\theoremstyle{definition}

\newtheorem{definition}[theorem]{Definition}

\theoremstyle{remark}

\newtheorem*{remark}{Remark}

\theoremstyle{remark}

\newtheorem*{remarks}{Remarks}

\begin{document}
\title{Singular Continuous Spectrum for the Laplacian on Certain Sparse Trees}
\author{Jonathan Breuer \\
\footnotesize Institute of Mathematics,  The Hebrew University of Jerusalem, \\
\footnotesize 91904 Jerusalem, Israel. \\
\footnotesize Email: jbreuer@math.huji.ac.il}
\date{}
\maketitle
\begin{abstract}
We present examples of rooted tree graphs for which the Laplacian has singular continuous spectral measures. 
For some of these examples we further establish fractional Hausdorff dimensions. 
The singular continuous components, in these models, have an interesting multiplicity structure. 
The results are obtained via a decomposition of the Laplacian into a direct sum of Jacobi matrices.
\end{abstract}

\section{Introduction}
This note deals with the spectral analysis of the discrete Laplacian on trees that have a certain 
sparseness in their coordination number (to be precisely defined below). We show that the 
spectral theory of the Laplacian on such trees bares similarities to the theory of 
one-dimensional Schr\"odinger operators with a sparse-barrier potential. In particular, this framework 
allows us to construct explicit examples of trees with singular continuous spectrum. Moreover, for some of these
models, the
spectral measures have fractional Hausdorff dimensions (see theorem \ref{growing-tree1} below). 
Graphs with singular continuous spectrum 
are known to exist \cite{graph-lap}; however, we are not aware of any previous explicit 
construction of a graph with non-trivial bounds on the spectral dimensions. What's more, we show 
that the singular continuous components occur with multiplicities that are related to the 
symmetry of the tree (theorem \ref{multiplicity}).

At the center of our analysis is a decomposition theorem (theorem \ref{decomposition}) for the 
Laplacian on a family of trees that exhibit a certain spherical symmetry. All our examples follow
from this decomposition by applying known methods from the theory of sparse-potential Schr\"odinger 
operators, mentioned above. Thus, for the applications, we content ourselves with giving the 
proper references and a few general remarks.

The paper is structured as follows. The basic result for sparse trees 
with singular continuous spectrum is described in section 2, along
with the decomposition theorem. 
The proofs of theorems \ref{decomposition} and \ref{multiplicity} are presented in Section 3. 
Further examples are given in section 4. 

\textbf{Acknowledgements}.
We are grateful to Michael Aizenman, Nir Avni, Yoram Last, Barry Simon, and Simone Warzel for useful discussions. 
We also wish to thank Michael Aizenman for the hospitality of Princeton where some of this 
work was done. 

This research was supported in part by THE ISRAEL SCIENCE FOUNDATION (grant no. 188/02) and by Grant no. \mbox{2002068} 
from the United States-Israel Binational Science Foundation (BSF), Jerusalem, Israel.

\section{Sparse Trees}
Recall that for a combinatorial tree, the distance between two vertices is defined as the number of 
edges of the unique path between them. 
\begin{definition}[Spherically Homogeneous Rooted Tree]
A rooted tree is called spherically homogeneous (SH) (see \cite{bass}) if any vertex $v$, at a distance $j$ from the root - 
$O$, is connected with 
$\kappa_j$ vertices at a distance $j+1$ from $O$. A locally finite (that is - the valence of 
every vertex is finite) 
spherically homogeneous tree is uniquely determined by the 
sequence $\{\kappa_j\}_{j=0}^\infty$. 
\end{definition} 

\setlength{\unitlength}{0.7cm}
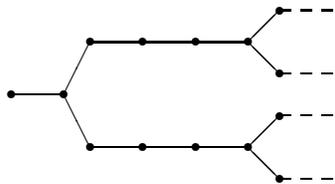
\begin{figure}[!hbp]
\begin{picture}(5,2)(0,-1)
\put(3,0){\line(1,0){1}}
\put(3,0){\circle*{0.15}}
\put(4,0){\circle*{0.15}}
\put(4,0){\line(1,2){0.5}}
\put(4.5,1){\circle*{0.15}}
\put(4,0){\line(1,-2){0.5}}
\put(4.5,-1){\circle*{0.15}}
\put(4.5,1){\line(1,0){1}}
\put(5.5,1){\circle*{0.15}}
\put(5.5,1){\line(1,0){1}}
\put(6.5,1){\circle*{0.15}}
\put(6.5,1){\line(1,0){1}}
\put(7.5,1){\circle*{0.15}}
\put(4.5,-1){\line(1,0){1}}
\put(5.5,-1){\circle*{0.15}}
\put(5.5,-1){\line(1,0){1}}
\put(6.5,-1){\circle*{0.15}}
\put(6.5,-1){\line(1,0){1}}
\put(7.5,-1){\circle*{0.15}}
\put(7.5,-1){\line(1,1){0.6}}
\put(7.5,-1){\line(1,-1){0.6}}
\put(7.5,1){\line(1,1){0.6}}
\put(7.5,1){\line(1,-1){0.6}}
\put(8.1,1.6){\line(1,0){0.2}}
\put(8.5,1.6){\line(1,0){0.2}}
\put(8.9,1.6){\line(1,0){0.2}}
\put(8.1,0.4){\line(1,0){0.2}}
\put(8.5,0.4){\line(1,0){0.2}}
\put(8.9,0.4){\line(1,0){0.2}}
\put(8.1,-1.6){\line(1,0){0.2}}
\put(8.5,-1.6){\line(1,0){0.2}}
\put(8.9,-1.6){\line(1,0){0.2}}
\put(8.1,-0.4){\line(1,0){0.2}}
\put(8.5,-0.4){\line(1,0){0.2}}
\put(8.9,-0.4){\line(1,0){0.2}}
\put(8.1,1.6){\circle*{0.15}}
\put(8.1,0.4){\circle*{0.15}}
\put(8.1,-1.6){\circle*{0.15}}
\put(8.1,-0.4){\circle*{0.15}}
\end{picture}
\caption{An example of a SH rooted tree with $\kappa_0=1$, $\kappa_1=2$, $\kappa_2=\kappa_3=\kappa_4=1$, $\kappa_5=2$ 
\ldots}
\end{figure}

Let $\{k_n\}_{n=1}^\infty$ be a sequence of natural numbers $>1$, 
and $\{L_n\}_{n=1}^\infty$ be a strictly increasing sequence of natural numbers. A spherically homogeneous rooted tree - 
$\Gamma$ - is said to be of type $\{L_n, k_n\}_{n=1}^\infty$, if 
\beq \label{kappa-j1}
\kappa_j= \left\{ \begin{array}{ll}
k_n & j=L_n \textrm{ for some } n\\
1 & \textrm{otherwise}
\end{array} \right.
\eeq
We shall say that $\Gamma$ is sparse if $(L_{n+1}-L_n) \rightarrow \infty$ rapidly, as $n \rightarrow \infty$.

Typical of the examples we construct is the following:
\begin{theorem} \label{k-bounded}
Let $k_0 \geq 2$ be a natural number and let $k_n \equiv k_0$. Assume that $(L_{n+1}-L_n) \rightarrow \infty$ and let 
$\Gamma$ be a SH rooted tree, of type $\{L_n, k_n\}_{n=1}^\infty$. 
Then the essential spectrum of $\Delta$ on $\Gamma$ contains the interval $[-2,2]$ and, provided 
$(L_{n+1}-L_n)$ increase sufficiently rapidly, any spectral measure for 
$\Delta$ is purely singular continuous on 
$(-2,2)$. By `sufficiently rapidly' we mean that $(L_{n+1}-L_n)$ has to be made sufficiently large with respect to 
$\{(L_{i+1}-L_i)\}_{i<n}$.
\end{theorem}

Note that since we are dealing with non-regular trees, there are two natural 
choices for the Laplacian:
\beq \label{lap1}
(\Delta f)(x)=\sum_{y : d(x,y)=1}f(y),
\eeq
and
\beq \label{lap2}
(\ti{\Delta} f)(x)=\sum_{y : d(x,y)=1 }f(y)-\#\{y : d(x,y)=1 \}\cdot f(x)
\eeq
where $\#A$, for a finite set $A$, is the number of elements in $A$ ($d(x,y)$ denotes the distance on the tree). 
Although we formulate all our results for $\Delta$, we note that they hold for $\ti{\Delta}$ as well, (with 
$(-2,2) \subseteq \bbR$ 
replaced by $(-4,0)$ where necessary). 

It is clear that if $\{k_n\}_{n=1}^\infty$ is a bounded sequence, then $\Delta$ 
on the tree is bounded and self-adjoint. For unbounded coordination number, the operator is unbounded and the issue 
of self-adjointess has to be addressed. 
\begin{definition} \label{normal-tree}
We call a SH rooted tree - $\Gamma$  - \emph{normal} if $\{k_n\}$ unbounded implies that 
$\limsup_{n \rightarrow \infty} (L_{n+1}-L_n)>1$. 
\end{definition}
Standard methods imply that the Laplacian on normal SH rooted trees is self-adjoint.

The following decomposition theorem allows us to represent $\Delta$ on 
$\Gamma$ as a direct sum of Jacobi matrices
\beq \label{jacobi}
J(\{a(j)\}_{j=1}^\infty,\{b(j)\}_{j=1}^\infty)=\left( \begin{array}{ccccc}
b(1)    & a(1) & 0      & 0      & \ldots \\
a(1)    & b(2) & a(2)    & 0      & \ldots \\
0      & a(2) & b(3)    & a(3)    & \ddots \\
\vdots & \ddots   & \ddots & \ddots & \ddots \\
\end{array} \right)
\eeq
with \beq \no  b(j) \in \bbR, \ a(j)>0. 
\eeq 
 
\begin{theorem} \label{decomposition}
Let $\Gamma$ be a normal rooted SH tree of type $\{L_n,k_n\}_{n=1}^\infty$.
Let 
\beq \label{M-n}
M_n= \left\{ \begin{array}{ll}
\prod_{j=1}^n k_j-\prod_{j=1}^{n-1}k_j & n>1\\
k_1-1 & n=1\\
1 & n=0.
\end{array} \right.
\eeq
Furthermore, let $R_0=0$ and $R_n=L_n+1,$ for $n\geq1$.
Then $\Delta$ is unitarily equivalent to a direct 
sum of Jacobi matrices, each operating on a copy of $\ell^2(\bbZ^+)$:
\beq \label{direct-sum-equivalence1}
\Delta \cong \oplus_{n=0}^\infty (\underbrace{J_n \oplus J_n \oplus \cdots \oplus J_n}_{M_n \textrm{ times}}) 
\eeq 
where $J_n=J(\{a_n(j)\}_{n=1}^\infty,\{b_n(j)\}_{n=1}^\infty)$ with 
\beq \label{a-n-j}
a_n(j)= \left\{ \begin{array}{ll}
\sqrt{k_m} & j=R_m-R_n \textrm{ for some } m>n\\
1 & \textrm{otherwise}
\end{array} \right.
\eeq
and
\beq \label{b-n-j}
b_n(j) \equiv 0.
\eeq
\end{theorem}
\begin{remarks}
1. For the case of a regular tree, a similar decomposition was discussed in \cite{froese, romanov} (see also \cite{solomyak} for a 
related result in the case of a metric tree).

2. As noted above, this theorem holds for $\ti{\Delta}$ as well, albeit with different values for $\ti{b}_n(j)$.
\end{remarks}

The next theorem is a simple corollary of theorem \ref{decomposition}. 
Its conditions are satisfied by all the cases we consider in this paper. 
\begin{theorem} \label{multiplicity}
Let $\Gamma$ be a normal SH rooted tree of type $\{L_n,k_n\}_{n=1}^\infty$. Consider $\Delta$ on $\Gamma$. 
Let $I \subseteq \bbR$ be an interval such that all the spectral measures restricted to $I$ are singular 
continuous. Let $P_I$ be the spectral projection onto $I$ (associated with 
$\Delta$), and let $M_n$ be as defined in \eqref{M-n}.
Then, $P_I \left( \ell^2(\Gamma) \right)$ decomposes as a direct sum of invariant spaces, 
$\oplus_{n=0}^\infty \mathcal{H'}_n$, 
such that $\Delta$, restricted to $\mathcal{H'}_n$, has uniform multiplicity $M_n$ and the measure classes 
associated with the representation of $\Delta$ restricted to $\mathcal{H'}_n$ are mutually 
disjoint.
\end{theorem}

\begin{remark}
It is not hard to show that the numbers $M_n$ are dimensions of certain irreducible 
representations of the symmetry group of $\Gamma$.
\end{remark}

Theorem \ref{decomposition} makes it clear that the spectral analysis of the Laplacian on sparse trees reduces to that of 
Jacobi matrices with `bumps' that are sparse along the subdiagonal and superdiagonal. Such matrices are analogous to 
discrete one-dimensional Schr\"odinger operators with 
potentials composed of sparse barriers. There is extensive literature on the spectral 
theory of such operators (see \cite{last-review} for a review of the relevant theory), 
showing that such potentials give rise to a 
variety of interesting spectral phenomena.
As noted in the introduction, all the examples we present in this paper are 
obtained by applying the (suitably modified) methods of the diagonal 
sparse case to the off diagonal case and using theorem \ref{decomposition}.
In particular, theorem \ref{k-bounded} follows from  
theorem \ref{decomposition} by the methods in \cite{pearson}.


\section{Decomposing the Laplacian}
We start with some terminology and notation: We use the shorthand $|v| \equiv d(v,O)$.
For any $v \in \scrV(\Gamma)$ we call \emph{the forward subtree of $v$} - $\Gamma_v$, the subtree of $\Gamma$ 
all
of whose vertices, $u$, satisfy the following two conditions 
\begin{enumerate}
\item $|v| \leq |u|$.
\item any vertex $v'$ on the unique path connecting $v$ and $u$ satisfies $|v| \leq |v'|$.
\end{enumerate}
We shall use $S_\Gamma(r) \equiv \{v \in \scrV(\Gamma) \mid |v| = r \}$.

\begin{proof}[Proof of Theorem \ref{decomposition}]
We shall decompose $\mathcal{H}=\ell^2(\Gamma)$ as a direct sum of spaces - 
$\oplus_{n=0}^\infty \mathcal{H}_n$, each 
invariant under $\Delta$, such that $\Delta$ restricted to 
$\mathcal{H}_n$ is unitarily equivalent to a direct sum of $M_n$ copies of $J_n$. We shall 
describe this decomposition inductively. 

We need a label for some of the vertices of $\Gamma$: At a distance $R_n$ from the root
 there are $\alpha_n \equiv \prod_{j=1}^n k_j$ vertices. These are naturally divided into 
$\alpha_{n-1}$ groups of $k_n$ vertices with common backward neighbor. 
We shall label the vertices on $S_\Gamma(R_n)$ by $\{v_{n,l}\}_{l=1}^{\alpha_n}$ where for each 
$m=1,2,\ldots,\alpha_{n-1}$, the vertices $\{v_{n,l}\}_{l=(m-1)k_n+1}^{mk_n}$ are all 
on the forward subtree of $v_{n-1,m}$.

In order to streamline the notation, we shall use $\delta_{n,l}$ for $\delta_{v_{n,l}}$ (the delta function at the 
vertex $v_{n,l}$). We shall also use $\Gamma_{n,l}$ for $\Gamma_{v_{n,l}}$. 

Now, let us define
$V_n=[\delta_{n,l}]_{l=1}^{\alpha_n}$ - the linear span of $\{\delta_{n,l}\}_{l=1}^{\alpha_n}$ -  
so that $\dim V_n=\alpha_n$. 
Let $\vphi_0=\delta_O$ and let 
\beq \label{H-0-def}
\mathcal{H}_0=\overline{[\Delta^n \vphi_0 \mid n=0,1,\ldots]}
\eeq
where $\overline{[\cdot]}$ for a linear subspace of $\mathcal{H}$ denotes its closure.  
An orthogonal basis for $\mathcal{H}_0$ is obtained by `Gram-Schmidting' the basis $\{\Delta^n \vphi_0\}$ - 
which results with normalized, radially symmetric functions supported on spheres around $O$ (the radial symmetry 
is a consequence of the spherical homogeneity of $\Gamma$). 
This implies that 
$\mathcal{H}_0$ is the subspace of radially symmetric functions. 

Now assume that we have defined $\mathcal{H}_n$ for $n=0,1,\ldots,(N-1)$ such that:

\emph{1.} $\mathcal{H}_i \perp \mathcal{H}_j$ for any $i \neq j$ and all spaces are invariant under $\Delta$ in the 

sense that $\Delta \left( D(\Delta) \cap \mathcal{H}_i \right) \subseteq \mathcal{H}_i$ 
(where $D(\Delta)$ is the domain of $\Delta$).

\emph{2.} For any vertex $v$ with $|v| < R_N$, $\delta_v \in \oplus_{n=1}^{(N-1)} \mathcal{H}_n$.

\emph{3.} For any $l>(N-1)$, $V_l \cap \left( \oplus_{n=1}^{(N-1)} \mathcal{H}_n \right)$ is $\alpha_{(N-1)}$ dimensional.

(All these properties hold trivially for $N-1=0$ if we let $\alpha_0=1$).
Recall that $M_n=\prod_{j=1}^n k_j-\prod_{j=1}^{n-1} k_j=\alpha_n-\alpha_{n-1}$. From the above, it follows that 
the orthogonal complement in $V_N$ (which is $\alpha_N$ dimensional) to $\oplus_{n=1}^{(N-1)} \mathcal{H}_n$ 
is $M_N$ dimensional and is spanned by 
$M_N$ mutually orthogonal unit vectors - $\vphi_{N,1},\ldots,\vphi_{N,M_N}$. Writing
\beq \label{induction-form1}
\vphi_{N,j}=\sum_{l=1}^{\alpha_N} a_{N,j}^l \delta_{N,l}, \qquad 1 \leq j \leq M_N
\eeq
and recalling that for all $m=1,2,\ldots,\alpha_{(N-1)}$, $\{v_{N,l}\}_{l=(m-1)k_N+1}^{mk_N}$ have a common backward 
neighbor on $S_\Gamma(L_N)$, we get (from \emph{2} above) that
\beq \label{basic-condition2}
\sum_{l=(m-1)k_N+1}^{mk_N} a_{N,j}^{l} = 0 
\eeq
for all $m$ and $j$.
Define
\beq \label{H-N-def}
\mathcal{H}_N=\overline{[\Delta^n \vphi_{N,j} \mid 1 \leq j \leq M_N \ n=0,1,\ldots]}. 
\eeq
Then we claim that
\beq \label{subspace-decomposition} 
\mathcal{H}_N=\oplus_{j=1}^{M_N} \mathcal{H}_{N,j} 
\eeq
with
\beq \label{H-Nj-def}
\mathcal{H}_{N,j}=\overline{[\Delta^n \vphi_{N,j} \mid n=0,1,\ldots]}.
\eeq
Indeed, \eqref{basic-condition2} together with the spherical homogeneity of $\Gamma$, implies that the orthogonal 
basis obtained from the Gram-Schmidt process applied to
$\{\Delta^n \vphi_{N,j}\}$,
is made of functions supported on $\{ v \mid |v| \geq R_N \}$ and having the form
\beq \label{GS-form2}
\frac{1}{\rho_{N,j}(r)}\sum_{|v|=r} a_{N,j}^v \delta_v
\eeq
where $r \geq R_N$, $\rho_{N,j}(r)>0$ is a normalizing factor and $a_{N,j}^v=a_{N,j}^l$ if $v \in \Gamma_{N,l}$. 
Together with the orthogonality of the various $\vphi_{N,j}$, this means that $\mathcal{H}_{N,j} \perp \mathcal{H}_{N,i}$ 
if $i \neq j$. Thus we see that properties \emph{1-3} above hold for $n=0,1,\ldots,N$. 

Having constructed $\mathcal{H}_n$ for all $n$ in this fashion, we get from properties \emph{1} and \emph{2} that indeed
\beq \label{H-decomposition}
\mathcal{H}=\oplus_{n=1}^\infty \mathcal{H}_n
\eeq 
and that each of the spaces in the direct sum is invariant under $\Delta$ in the sense that 
$\Delta \left( D(\Delta) \cap \mathcal{H}_n \right) \subseteq \mathcal{H}_n$. This almost means that $\Delta$ decomposes 
as a direct sum. The only missing point is that if $\Delta$ is unbounded, the above does not necessarily mean that 
$\mathcal{H}_n$ is invariant under $(\Delta-z)^{-1}$ for $z \in \bbC \setminus \bbR$. However, it is easy to see that the 
moment problem associated with the operation of $\Delta$ on $\vphi_{n,j}$ (any $n$ and $1 \leq j \leq M_n$) is determinate, 
so it follows from proposition 4.15 in \cite{moment} that $\mathcal{H}_n$ is 
indeed invariant in both senses. Thus, we have that
\beq \label{Delta-decomposition}
\Delta=\oplus_{n=1}^\infty \Delta_n
\eeq
with $\Delta_n$ denoting the corresponding restricted operators.

In order to complete the proof, we need to show that 
\beq \label{Delta-n}
\Delta_n \cong \underbrace{J_n \oplus J_n \oplus \cdots \oplus J_n}_{M_n \textrm{ times}}.
\eeq
Knowing \eqref{subspace-decomposition} and \eqref{GS-form2}, however, this is now a matter of 
simple computation 
(since it is easy to see that 
$\rho_{n,j}(r)=\sqrt {\# \left(S_\Gamma(r) \cap \Gamma_{n,l} \right)}$ for any $1 \leq l \leq \alpha_n$ and $r \geq R_n$).

\end{proof}

\begin{proof}[Proof of Theorem \ref{multiplicity}]
Theorem \ref{decomposition} says that it suffices to consider the measures $\mu_n$ - the spectral measures of $J_n$ and 
$\delta_1 \in \ell^2(\bbZ^+)$ (since this is a cyclic vector). Obviously, these measures occur with multiplicity at least 
$M_n$, so we only need to show that their singular continuous parts are mutually singular. Note now that for $n_1>n_2$, 
$J_{n_1}$ can 
be obtained from $J_{n_2}$ by `stripping off' the $(R_{n_1}-R_{n_2})$ leftmost columns and the same number of rows from the 
top. The fact that the singular continuous part of $\mu_{n_1}$ is singular, with respect to the singular continuous part of 
$\mu_{n_2}$, follows, now, from the characterization of the appropriate supports in terms of $m$-functions 
(see e.g.\ \cite{rankone}) and from the continued fraction expansion of $m$ \cite{m-function}.
(The spaces $\mathcal{H'}_n$ are just $P_I(\mathcal{H}_n)$).
\end{proof}

\section{Singular Continuous Spectrum for Sparse Trees}  

Let $\Gamma$ be a sparse, normal SH rooted tree. 
As noted in the proof of theorem \ref{multiplicity}, all the matrices $J_n$, in the decomposition of the Laplacian, are 
actually various `tails' of $J_0$. Since, in the sparse models, it is the asymptotics that determine the spectral type, this 
means that the spectral analysis of the Laplacian reduces essentially to the analysis of a single Jacobi matrix.

A good reason for considering trees with unbounded $\{k_n\}$, is the fact that for trees with $k_n \rightarrow \infty$, it 
is easy to state explicit growth conditions on $\{L_n\}$ which make the spectrum singular continuous. In particular, 
$k_n \rightarrow \infty$ implies absence of absolutely continuous spectrum \cite{last-simon}, 
so the following is a straightforward adaptation of an idea of Simon-Stolz \cite{simon-stolz} 
to our case:

\begin{theorem} \label{k-unbounded}
Let $\{k_n\}_{n=1}^\infty$ be a sequence of natural numbers such that $k_n \rightarrow \infty$ as $n \rightarrow \infty$.
Let $\alpha_n=\prod_{j=1}^n k_j$. 
Assume that $(L_{n+1}-L_n) \rightarrow \infty$ and let $\Gamma$ be a SH rooted tree of type $\{L_n,k_n\}_{n=1}^\infty$.
Then the spectrum of $\Delta$ on $\Gamma$ consists of the interval $[-2,2]$ along with some discrete point spectrum 
outside this interval. If for some $\varepsilon>0$, 
\beq \label{sc-condition}
\limsup_{n \rightarrow \infty} \frac{(L_{n+1}-L_n)}{\alpha_n^{(1+\varepsilon)}}>0,
\eeq 
then any spectral measure for $\Delta$ is purely singular continuous on $(-2,2)$. 
\end{theorem}

As certain sparse potentials have been constructed with spectral measures of fractional Hausdorff 
dimensionality, it seems natural to try to construct trees with this property as well. An adaptation of an example of 
Jitomirskaya-Last \cite{jit-last} (see also \cite{tcherem}) achieves just that:

\begin{theorem} \label{growing-tree1}
Let $L_n=2^{(n^n)}$.
Let $\beta>0$ and $k_n=\left[ L_n^\beta \right]$. Let $\Gamma_\beta$ be the corresponding tree. 
Then the restriction to $(-2,2)$ of any spectral measure for $\Delta$ on $\Gamma_\beta$, is 
supported on a set of Hausdorff dimension $\frac{2}{2+\beta}$ and does not give weight to sets
of Hausdorff dimension less than $\frac{1}{1+\beta}$.

Letting $k_n=\left[ L_n^{\beta_n(c)} \right]$ with $\beta_n(c)=c\frac{(n+1)^{(n+1)}}{n^n}$ for some $1>c>0$, 
we get that any spectral measure on $(-2,2)$ is purely singular continuous and supported on a 
set of Hausdorff dimension $0$. 
\end{theorem}


\end{document}